\documentclass[11pt]{article}
\title{The Rationality of Sol--Manifolds}
\author{Thomas A. Putman}
\usepackage{amsmath}
\usepackage{amsfonts} 
\usepackage{amssymb}

\newtheorem{theorem}{Theorem}[section]
\newtheorem{lemma}[theorem]{Lemma}
\newtheorem{corollary}[theorem]{Corollary}

\newenvironment{proofof}[1]{{\bf {\smallskip}{\noindent}Proof of #1: }}{$\square$ \smallskip}
\newenvironment{proof}{{\bf {\smallskip}{\noindent}Proof: }}{$\square$ \smallskip}
\newenvironment{definition}{{\bf {\smallskip}{\noindent}Definition: }}{\smallskip}
\newenvironment{remark}{{\bf {\smallskip}{\noindent}Remark: }}{\smallskip}
\newenvironment{acknowledgements}{{\bf {\smallskip}{\noindent}Acknowledgments: }}{\smallskip}
\newenvironment{history}{{\bf {\smallskip}{\noindent}History and Comments: }}{\smallskip}
\newenvironment{example}{{\bf {\smallskip}{\noindent}Example: }}{\smallskip}
\newenvironment{claim}[1]{{\bf {\smallskip}{\noindent}Claim #1: }}{}
\newenvironment{claimproof}{{\bf {\noindent}Proof of Claim: }}{$\square$}

\newcommand\HBolic{\mbox{$\mathbb{H}$}}
\newcommand\SL{\mbox{SL}}
\newcommand\PSL{\mbox{PSL}}
\newcommand\R{\mbox{$\mathbb{R}$}}
\newcommand\Z{\mbox{$\mathbb{Z}$}}
\newcommand\Max{\mbox{max}}
\newcommand\Min{\mbox{min}}
\newcommand\Th{\mbox{th}}
\newcommand\type{\mbox{type}}
\newcommand\unreducedtype{\mbox{$\widetilde{\text{type}}$}}
\newcommand\reducedtype{\mbox{type}}
\newcommand\Trace{\mbox{trace}}
\newcommand\Growth{\mbox{${\mathcal G}$}}
\newcommand\Center{\mbox{${\mathcal C}$}}
\newcommand\Tail{\mbox{${\mathcal T}$}}
\newcommand\TailLen{\mbox{$\overline{\mathcal T}$}}
\newcommand\Head{\mbox{${\mathcal H}$}}
\newcommand\HeadLen{\mbox{$\overline{\mathcal H}$}}
\newcommand\Rev{\mbox{rev}}
\newcommand\Abelian{\mbox{ab}}

\begin{document}

\maketitle

\begin{abstract}
Let $\Gamma$ be the fundamental group of a manifold modeled on 3-dimensional
Sol geometry.  We prove that $\Gamma$ has a finite index subgroup $G$ which
has a rational growth series with respect to a natural generating set.  We do
this by enumerating $G$ by a regular language.  However, in contrast to most
earlier proofs of this sort our regular language is not a language of words
in the generating set, but rather reflects a different geometric structure
in $G$.
\end{abstract}

\section{Introduction}
\label{section:introduction}

Let $\Gamma$ be a group with a finite generating set $S$.  For $g \in \Gamma$, 
let $\|g\|$ be equal to the length of the shortest word in $S \cup S^{-1}$
representing $g$, and for $g_1,g_2 \in \Gamma$ set 
$d(g_1,g_2)=\|g_1^{-1}g_2\|$.  This is known as the {\em word metric\/} on 
$\Gamma$.  The growth of the size of balls in this metric constitutes a 
central object of study in geometric group theory (see \cite[Chapters 6--7]{Ha}
for a survey).

To study the growth of $\Gamma$, it is natural to define the {\em growth 
series\/} of $\Gamma$ to be the power series
$$\Growth(\Gamma)=\sum\nolimits_{i=0}^{\infty}c_i z^i$$
where $c_i=|\{g \in \Gamma:\|g\|=i\}|$.  In many cases, it turns out that 
$\Growth(\Gamma)$ is a rational function.  The first nontrivial example of this is 
in \cite{Bo}, where an exercise outlines a proof that all Coxeter groups have 
rational growth with respect to a Coxeter generating set.  Perhaps the most 
remarkable theorem of this type is in Cannon's paper \cite{Ca}, which proves 
that all word hyperbolic groups have rational growth with respect to 
{\em any\/} finite generating set (\cite{Ca} only proves this for fundamental groups 
of compact hyperbolic manifolds, but it contains all the ideas necessary for 
the extension to word hyperbolic groups -- see \cite{CDP} for a complete 
account).  

In this paper we study the growth series of the fundamental groups $\Gamma$ of 
torus bundles over the circle with Anosov monodromy.  In other words, 
$\Gamma=\Z^2 \rtimes_M \Z$ with $M \in \SL_2(\Z)$ a matrix with two distinct 
real eigenvalues.  These are the fundamental groups of 3--manifolds modeled on 
Sol geometry.  Our main theorem is the following:

\begin{theorem}[Main Theorem]
\label{theorem:maintheorem}
Let $\Gamma$ be the fundamental group of a 3--dimensional Sol manifold.  Then 
there exists a finite index subgroup $G$ generated by a finite set $S$ so that 
$G$ has rational growth with respect to $S$.  In other words, $\Gamma$ is 
virtually rational.
\end{theorem}

This theorem is part of two different streams of research.  On the one hand, 
there have been many papers investigating the growth series for lattices in 
Thurston's eight 3--dimensional model geometries (see \cite{Be1, Be2, Ca, NS1, 
Sh1, Sh2, We}).  After Theorem \ref{theorem:maintheorem}, the only remaining 
geometry for which there is not some general theorem is $\widetilde{\SL_2}$, 
although some progress has been made on this case by Shapiro \cite{Sh2}.

On the other hand, there has also been significant research on the growth
series of finitely generated solvable groups.  Kharlampovich has produced
a 3--step solvable group which has an unsolvable word problem \cite{Kh}.
Since all groups with rational growth series have a solvable word problem
(the rational growth series allows one to calculate the size of balls in
the Cayley graph, which one can then construct using a brute force
enumeration), it follows that Kharlampovich's example does not have rational
growth with respect to any set of generators.

One can therefore hope for general results only for 1 and 2--step solvable 
groups.  The 1--step solvable groups are the finitely generated abelian groups.
Benson has proven that more generally all finitely generated virtually abelian 
groups have rational growth with respect to any finite set of generators \cite{Be1}.  
The 2--step solvable groups are divided into the nilpotent and non-nilpotent 
groups.  A fundamental set of examples of 2--step nilpotent groups are the 
lattices in 3--dimensional Nil geometry.  These correspond to groups of the 
form $\Z^2 \rtimes_M \Z$ with $M \in \SL_2(\Z)$ a matrix with two different 
non-real eigenvalues, which necessarily must lie on the unit circle.  Benson, 
Shapiro, and Weber \cite{Be2, Sh1, We} have shown that lattices in 
3--dimensional Nil geometry have rational growth with respect to a certain 
generating set.  This has been generalized by Stoll \cite{Sto}, who showed that
all 2--step nilpotent groups with infinite cyclic derived subgroup have 
rational growth with respect to some generating set.  He also showed that many 
such groups (those with ``Heisenberg rank at least $2$'') have transcendental 
growth with respect to some other generating set.  This demonstrates that, in 
contrast to most natural group properties studied by geometric group theorists,
rational growth can depend strongly on the choice of generating set.

The non-nilpotent solvable case is divided into the polycyclic and the 
non-polycyclic cases.  A fundamental set of examples of 2--step solvable 
non-polycyclic groups are the solvable Baumslag--Solitar groups $BS(1,n)$.  
Brazil and Collins--Edjvet--Gill have shown that these have rational growth 
with respect to the standard set of generators \cite{Bz,CEG}.  A fundamental 
set of examples of 2--step solvable polycyclic groups are the torsion--free 
abelian--by--cyclic groups.  These are groups of the form $\Z^n \rtimes_M \Z$ 
with $M \in \SL_n(\Z)$.  Theorem \ref{theorem:maintheorem} combined with the 
result of Benson, Shapiro, and Weber referred to in the previous paragraph 
cover the case n=2, with ours being the ``generic'' case since our eigenvalues 
do {\em not\/} lie on the unit circle.

The strategy of our proof is as follows.  The subgroup $G$ we consider is generated
by two elements $a$ and $t$.  We first define a surjective function which 
associates to a word $w$ in the free group on $a$ and $t$ a pair $(\unreducedtype(w),h(w))$  
where $\unreducedtype(w) \in \Z[z,z^{-1}]$ and $h(w) \in \Z$.  The word metric
on the free group induces a ``size function'' on such pairs for which
there is a rather simple formula.  It turns out that two words $w_1$ and
$w_2$ determine the same element of $G$ if and only if $h(w_1)=h(w_2)$ and
$\unreducedtype(w_1)=\unreducedtype(w_2)$ modulo a certain principal ideal $I$ of
$\Z[z,z^{-1}]$.  It is easy to construct a ``size-preserving'' enumeration of 
$\Z[z,z^{-1}] \oplus \Z$ by a regular language $L$.  The map
$$\Z[z,z^{-1}] \oplus \Z \longrightarrow (\Z[z,z^{-1}]/I) \oplus \Z$$
induces a ``quotient'' $L/P$ of $L$ with an induced ``size function''.  We
conclude by proving that $L/P$ satisfies a certain negative curvature-like condition 
(the {\em falsification by fellow traveler property}).  This allows us to
enumerate $L/P$ by a regular language, which by well-known results is
enough to prove that it (and therefore $G$) has a rational growth series.

\begin{remark}
Though in theory our methods are entirely constructive, in practice the
finite state automata we build are so huge that it is impractical to
calculate any examples.
\end{remark}

\begin{history}
In his unpublished thesis \cite{Gr}, Grayson claimed to prove Theorem
\ref{theorem:maintheorem} whenever the trace of the monodromy is even.  
However, his proof is insufficient (see the remarks in Section 
\ref{section:typesandheights} for a more detailed discussion).  Our methods are 
rather different from his methods.  He attempts to write down a complicated 
recurrence relation between balls of different radii.  As indicated above, we 
instead use the theory of finite state automata.  In addition, the generating
sets we use are slightly different from his generating sets.  We do, however, use some of 
his ideas.  In particular, he introduced the notions of {\em types\/} and 
{\em heights\/} described in Section \ref{section:typesandheights} (though he did 
not distinguish between the reduced and unreduced types), and the elegant proof
of Theorem \ref{theorem:graysontheorem} is due to him.

After this paper was complete, we learned that in an unpublished paper Parry 
had given a proof of Theorem \ref{theorem:maintheorem}, following Grayson's 
basic outline \cite{P}.  Like Grayson, he assumes that the trace of the 
monodromy is even.  However, Parry was able to use a computer to calculate 
some growth functions explicitly.  We reproduce the result of his calculation 
in Section \ref{section:questions}.
\end{history}

\begin{acknowledgements}
I would like to thank my advisor Benson Farb for introducing me to this problem 
and for offering many corrections to previous versions of this paper.  I would
like to thank Chris Hruska for several useful conversations and for 
helping me to make sense of Grayson's work.  In addition, I would like to thank
Murray Elder for providing some useful references and offering some 
corrections and Walter Parry for providing me with his unpublished manuscript \cite{P}.
Finally, I wish to thank an anonymous referee for offering many detailed corrections
and helping me greatly improve the exposition.
\end{acknowledgements}

\subsection{Outline and Conventions}
In Section \ref{section:preliminaries}, we review some preliminary material on 
Sol manifolds, regular languages, etc.  Next, in Section 
\ref{section:partitioning} we discuss a technical condition on partitions of 
regular languages which implies rational growth.  This condition, the 
{\em falsification by fellow traveler\/} property, is inspired by but different
from the condition of the same name defined by Neumann--Shapiro in \cite{NS1}.
Section \ref{section:typesandheights} is then devoted to the bijection
$$G \longrightarrow (\Z[z,z^{-1}]/I) \oplus \Z$$
discussed above.  Finally, in Section
\ref{section:language} we construct a sequence of regular languages $L_n$ which enumerate
$\Z[z,z^{-1}] \oplus \Z$ in a ``size-preserving'' manner, and we prove that the
partition $P_n$ of $L_n$ induced by the natural map
$$\Z[z,z^{-1}] \oplus \Z \longrightarrow (\Z[z,z^{-1}]/I) \oplus \Z$$
satisfies the falsification by fellow traveler property for sufficiently large $n$.  This
proves Theorem \ref{theorem:maintheorem}.  We conclude by discussing some open questions 
in Section \ref{section:questions}.

We will frequently manipulate {\em Laurent polynomials} over $\Z$, that is
elements of $\Z[z,z^{-1}]$.  When we refer to such a polynomial as $\sum\nolimits_j c_j z^j$,
we mean that all but finitely many of the $c_j$ equal $0$.

\section{Preliminaries}
\label{section:preliminaries}
\subsection{Sol Manifolds}
As discussed in \cite{Th}, 3-dimensional Sol manifolds are 2-dimensional torus bundles 
over the circle whose monodromy $M \in \SL_2(\Z)$ is {\em Anosov}, that is $M$ has two 
distinct real eigenvalues.  Equivalently, $|\Trace(M)|>2$.  Let $M$ be such a 
matrix, and let $a$ and $b$ be the standard generators for $\Z^2$.  Hence, $Ma$ and $Mb$ are
well defined.  Abusing notation in the obvious way, we say that the {\em torus 
bundle group\/} with monodromy $M$ is the group with the presentation
$$\Gamma = \langle a, b, t \, | \, [a,b]=1, tat^{-1}=Ma, tbt^{-1}=Mb \rangle$$
Observe that $G=\langle a,t \rangle$ is a finite index subgroup of $\Gamma$.  
Now, the minimal polynomial of the matrix $M$ is equal to $1-\Trace(M) z+z^2$.  Hence
$M^k a$ is in the lattice generated by $a$ and $M a$.  In other words, the group $G$ 
corresponds to the 2-dimensional torus bundle whose fiber is generated by $a$ and $M a$.
It is easy to see that $G$ is isomorphic to the torus bundle group with monodromy
$$\begin{pmatrix}
         0 & -1 \\
         1 & \Trace(M) \\
    \end{pmatrix}$$
We will prove that $G$ has rational growth with respect to the generating set 
$\{a,t\}$.

\subsection{Sized Sets and Languages}
\label{section:sizedsets}
In the course of our proof, we will construct a series of objects whose growth 
reflects the growth series of $G$.  The ``size functions'' on these objects 
come from very different sources.  The following formalism provides a language 
with which to compare these objects:

\begin{definition}
A {\em sized set\/} is a set $X$ together with a size function 
$\|\cdot\|:X \longrightarrow {\Z}_{\geq 0}$. 
\end{definition}

\noindent 
Set $c_i=|\{x \in X : \|x\|=i\}|$.  We will only consider
sized sets with $c_i < \infty$ for all $i$.  There is therefore
an associated generating function 
$$\Growth(X)=\sum\nolimits_{i=0}^{\infty} c_i z^i$$

\begin{definition}
Let $X$ be a sized set and $P$ be a partition of $X$.  In other words, $P$
is a set of pairwise disjoint subsets of $X$ so that
$$\bigcup_{A \in P}A=X$$
We define $X/P$ to be the sized set whose elements are elements of $P$ and 
whose size function is
$$\|A\|=\text{min}\{\|x\| : x \in A\}$$
If $x \in X$, then we will denote by $\overline{x}$ the set $A \in X/P$  
with $x \in A$.  If $x,y \in X$ satisfy $\overline{x}=\overline{y}$, then
we will say that $x$ equals $y$ modulo $P$.
\end{definition}

\begin{definition}
Let $(X_1,\|\cdot\|_1)$ and $(X_2,\|\cdot\|_2)$ be sized sets.  A bijection
$\psi : X_1 \longrightarrow X_2$ is a {\em near-isometry\/} if there
is some constant $c$ so that for all $x \in X_1$ we have 
$\|x\|_1=\|\psi(x)\|_2+c$.  If $c=0$, then a near-isometry is an {\em isometry}.
\end{definition}

\begin{remark}
Observe that if $X_1$ and $X_2$ are near-isometric with respect to
a constant $c$ then $\Growth(X_1) = z^c \Growth(X_2)$.  In particular,
$\Growth(X_1)$ is a rational function if and only if $\Growth(X_2)$ is.  Our
use of near-isometries is purely a matter of convenience -- they allow
us to have a somewhat simpler definition of the languages $L_n$ we construct
in Section \ref{section:language}.
\end{remark}

\noindent
Our primary source of sized sets will be the following:

\begin{definition}
Let $A$ be a finite set, which we will call the {\em alphabet}.  A 
{\em language\/} $L$ over $A$ is a subset of $A^\ast$, the set of finite 
sequences of elements of $A$.  Elements of $L$ are called {\em words}.
\end{definition}

\noindent
Languages can be considered sized sets in the following way.  Let $L$ be a 
language over $A$.  Consider some $\phi:A \longrightarrow \Z_{> 0}$, which 
we will call the {\em weighting}.  For $a_1 a_2 \cdots a_k \in L$, define
$$\|a_1 a_2 \cdots a_k\|=\sum\nolimits_{i=1}^k \phi(a_i)$$

\begin{example}
Let $H$ be a group with a finite set of generators $S$.  Let $L$ be the 
language of all words in $S \cup S^{-1}$ with weighting $1$ for each generator.  Finally, 
let $P$ be the partition which identifies two words if they represent the same 
element in $H$.  The series $\Growth(L/P)$ is then the usual growth series for $H$.
\end{example}

\subsection{Regular Languages}
\label{section:regularlanguages}
We quickly review the theory of finite state automata and regular languages.  
For more details see, e.g., \cite[Chapter 1]{ECHLPT}.

\begin{definition}
A {\em finite state automaton\/} on $n$ strings is a 5--tuple
$$(A,(V,E),S,F,l)$$
with $A$ a finite set (called the {\em alphabet}), $(V,E)$ a finite directed 
graph (called the {\em state graph}), $S \in V$ (called the {\em start state}), 
$F \subset V$ (called the {\em final states}), and  
$l:E \longrightarrow \prod\nolimits_{i=1}^n (A \cup \{\$\})$ with $\$$ some symbol 
disjoint from $A$ ($l$ is called the {\em transition label}; ``$\$$'' is a 
symbol for the end of a word) satisfying the following condition : if 
$l(e)=(\ldots,\$,\ldots)$ with the $\$$ in the $k^{\Th}$ place, then $l(f)$ 
also has a $\$$ in the $k^{\Th}$ place for all edges $f$ so that there is a 
finite (oriented) path
$$e=e_0,e_1,\ldots,e_m=f$$
\end{definition}

\begin{definition}
Let $Z=(A,(V,E),S,F,l)$ be a finite state automaton on $n$ strings.  We define
the language $L(Z) \subset \prod\nolimits_{i=1}^n A^{\ast}$ to be the following.  
Consider any element $(w_1,\ldots,w_n) \in \prod\nolimits_{i=1}^n A^{\ast}$.  Assume 
that the longest word in this tuple has $m$ letters.  For $1 \leq i \leq n$ and 
$1 \leq j \leq m$ define $w_i^j$ to be the $j^{\Th}$ letter of $w_i$ if $j$ is 
at most the length of $w_i$ and $\$$ otherwise.  Then 
$(w_1,\ldots,w_n) \in L(Z)$ if and only if there is some path
$$S=v_1,e_1,v_2,e_2,\ldots,e_{m},v_{m+1} \in F$$
so that for $1 \leq j \leq m$ we have
$$l(e_j) = (w_1^j,w_2^j,\ldots,w_n^j)$$
We say that $L(Z)$ is a {\em regular language}.
\end{definition}

\begin{remark}
Observe that we are abusing the word ``language'' in this definition : only the
case $n=1$ is an actual language.  
\end{remark}

\begin{remark}
One should think of this as a machine able to keep track of a finite amount of 
information.  The vertices of the state graph correspond to the different
states in which the machine can be, and the machine moves from the state $s_1$
to the state $s_2$ upon reading $\alpha$ if there is an edge $e$ between $s_1$
and $s_2$ with $l(e)=\alpha$.
\end{remark}

\noindent
The following theorem demonstrates the flexibility of regular languages:

\begin{theorem}{\cite[Proposition 1.1.4, Theorem 1.2.8, Corollary 1.4.7]{ECHLPT}}
\label{theorem:regular}
The class of regular languages is closed under all first order predicates 
(i.e. $\cup$, $\cap$, $\neg$, $\forall$, and $\exists$) and under concatenation.  In
addition, if $L$ is a regular language on $n$ strings then the following language is regular:
\begin{align*}
\Rev(L)=\{&(a_{1,1} \cdots a_{1,m_1},\ldots,a_{n,1} \cdots a_{n,m_n}) : \\
          &(a_{1,m_1} \cdots a_{1,1},\ldots,a_{n,m_n} \cdots a_{n,1}) \in L\}
\end{align*}
\end{theorem}

\noindent 
We will also need the following theorem:

\begin{theorem}
\label{theorem:regularrational}
Let $Z=(A,(V,E),S,F,l)$ be a finite state automaton on one string and let 
$\phi : A \longrightarrow \Z_{> 0}$ be a weighting.  Then the generating 
function $\Growth(L(Z))$ with the language $L(Z)$ weighted by $\phi$ is a rational 
function.
\end{theorem}
\begin{proof}
It is a standard fact (see, for instance, \cite[Theorem 9.1]{CDP}) that 
$\Growth(L(Z))$ is rational if $\phi$ is the constant function $1$.  To deduce 
the general case from this, replace each edge $e$ in $(V,E)$ by a path of 
length $\phi(l(e))$ with each edge in the path labeled by $l(e)$.
\end{proof}

\begin{remark}
When we refer to a regular language without specifying how many strings it has,
we are referring to a regular language on one string.
\end{remark}

\section{Partitioning Regular Languages}
\label{section:partitioning}
Fix a regular language $L$ with weighting $\phi$.  Consider a partition $P$ 
of $L$.  By Theorem \ref{theorem:regularrational}, we know that $L$ has
a rational generating function.  In this section we give a sufficient condition
for $L/P$ (see Section \ref{section:sizedsets} for the definition of $L/P$) to 
have a rational generating function.  Our condition, the {\em falsification by 
fellow traveler\/} property, allows us to construct a regular sublanguage of 
$L$ containing exactly one word of minimal size from each set in $P$.  It is inspired
by the property of the same name in \cite{NS1}.  We begin
with two preliminary definitions. 

\begin{definition}
We say that $L/P$ has a {\em regular cross section\/} if there is some regular
sublanguage $L' \subset L$ so that for all $A \in P$ there is a unique
$x \in L'$ with $x \in A$.  If in addition all such $x$ satisfy
$$\|x\| = \Min\{\|x'\| : x' \in A\}$$
then we say that $L'$ is a {\em regular minimal cross section\/} of $L/P$.
\end{definition}

\begin{definition}
We say that a regular language $R \subset L \times L$ is an {\em acceptor} for
a partition $P$ of $L$ if  
$$(w,w') \in R \Longrightarrow \textit{ $\overline{w}=\overline{w}'$ and $(w',w) \in R$}$$
\end{definition}

\noindent
Our condition is the following:

\begin{definition}
We say that a partition $P$ with an acceptor $R$ 
has the {\em falsification by fellow traveler\/} property if 
there is some constant $K$ and some regular sublanguage $L'$ of $L$ containing 
at least one minimal size representative of each set in $P$ so that if 
$w \in L'$ is not a minimal size representative in $L/P$ then there is some word 
$w' \in L$ so that the following are true.
\begin{itemize}
\item $(w,w') \in R$ (and, in particular, $\overline{w}=\overline{w}'$)
\item $\|w'\| < \|w\|$
\item For any $j$, let $s$ and $s'$ be the initial segments of $w$ and $w'$ of 
length $j$.  Then $| \|s\| - \|s'\| | \leq K$ (the words $w$ and $w'$ are said to 
{\em $K$--fellow travel\/}).
\end{itemize}
We also require that if $w,w' \in L'$ are both minimal size representatives of the 
same element of $L/P$ then $(w,w') \in R$.
\end{definition}

\noindent
Our main theorem about such partitions is the following:

\begin{theorem}
\label{theorem:fellowtraveller}
Let $P$ be a partition of a weighted regular language $L$ with an acceptor $R$.
Assume that $P$ has the falsification by fellow
traveler property.  Then $L/P$ has a regular minimal cross section.
\end{theorem}

\noindent
Theorems \ref{theorem:fellowtraveller} and \ref{theorem:regularrational} 
imply the following:

\begin{corollary}
\label{corollary:fellowtravellerrational}
Let $P$ be a partition of a weighted regular language $L$ with an acceptor $R$ so that $P$ has
the falsification by fellow traveler property.  Then $L/P$ has rational growth.
\end{corollary}

\noindent
Before proving Theorem \ref{theorem:fellowtraveller}, we need a lemma.

\begin{lemma}
\label{lemma:fellowtraveller}
Let $P$ be a partition of a regular language $L$ with a weighting $\phi$ and an 
acceptor $R$, and let $K$ be a natural number.  Then the following 
language is regular:
\begin{align*}
L_K=\{ & (w_1,w_2) \in L \times L : (w_1,w_2) \in R, \|w_1\| > \|w_2\|, \\ 
      & \textit{and } w_1 \textit{ and } w_2 \textit{ $K$--fellow travel} \}
\end{align*}
\end{lemma}
\begin{proof}
Observe that $L_K$ is the intersection of $R$ and the language
\begin{align*}
L_K'=\{ & (w_1,w_2) \in (A^{\ast})^2 : \|w_1\| > \|w_2\| 
         \textit{ and } w_1 \textit{ and } \\
       & w_2 \textit{ $K$--fellow travel }\}
\end{align*}
By Theorem \ref{theorem:regular} it is therefore enough to show that $L_K'$ is 
regular.  We construct an automaton accepting $L_K'$ as follows.  For simplicity,
we will extend $\phi$ to $A \cup \{\$\}$ by setting $\phi(\$)=0$.  Our automaton
has $2K+1$ states labeled $-K,\ldots,K$ plus a failure state.  The label
on a numbered state represents the difference between the portions of $w_1$ and $w_2$
read thus far.  We begin in state $0$.  Now assume that we are in a state
$i$ and read $a$ from $w_1$ and $b$ from $w_2$.  If $|i+\phi(a)-\phi(b)|>K$,
then $w_1$ and $w_2$ have ceased to $K$-fellow travel, so we go to the failure state.  
Otherwise, we go to the state $i+\phi(a)-\phi(b)$.
We succeed and accept $(w_1,w_2)$ if we end in a state with a positive label, and 
we fail otherwise.
\end{proof}

\noindent
We now prove Theorem \ref{theorem:fellowtraveller}.

\begin{proofof}{\ref{theorem:fellowtraveller}}
Let $L'$ be the regular sublanguage of $L$ and $K$ be the constant given by the 
definition of the falsification by fellow traveler property.  By Lemma 
\ref{lemma:fellowtraveller}
\begin{align*}
L_K=\{ & (w_1,w_2) \in L \times L : (w_1,w_2) \in R, \|w_1\| > \|w_2\|, \\
       & \textit{and } w_1 \textit{ and } w_2 \textit{ $K$--fellow travel} \}
\end{align*}
is a regular language.  Hence by Theorem \ref{theorem:regular}
\begin{align*}
L''=\{ & w \in L' :\textit{ there does not exist any } w' \in L 
         \textit{ so that } \\
       & (w,w') \in L_K\}
\end{align*}
is a regular language.  This language is composed of minimal size representatives 
in $L/P$.  It contains at least one representative of each element.  By a 
remark on page 57 of \cite{ECHLPT}, the language
$$S=\{(w_1,w_2) \in L'' \times L'' : \textit{$w_1$ is short-lex less than $w_2$}\}$$
is regular (see \cite[p. 56]{ECHLPT} for the definition of the short-lex 
ordering.  For our purposes its only important property is that it is a total 
ordering on the set of words).  
We conclude from Theorem \ref{theorem:regular} that
$$L'''=\{w \in L'' : \textit{ for all $w' \in L''$ we have $(w',w) \notin S \cap R$}\}$$ 
is regular.  By the definition of the falsification by fellow traveler 
property, $L'''$ contains a {\em unique\/} representative of minimal length 
for each element of $L/P$; i.e. it is a regular minimal cross section of $L/P$.
\end{proofof}

\section{Types and Heights}
\label{section:typesandheights}
Fix a torus bundle group $\Gamma$ with monodromy $M$.  Recall that we are examining 
the finite index subgroup $G=\langle a,t \rangle$.

\subsection{Definitions}
Consider some $g \in G$.  Since $G \subset \Gamma = \Z^2 \rtimes_M \Z$, we
can regard $g$ as a pair $(x,h)$ with $x \in \Z^2$ and $h \in \Z$.  We
will call $h$ the {\em height} of $g$ (denoted $h(g)$) and $x$ the {\em type}
of $g$ (denoted $\type(g)$).  

Denote by $F_T$ the free group on a set $T$.  Consider $w \in F_{\{a,t\}}$ which maps
to $\overline{w} \in G$.  Set $h(w)=h(\overline{w})$ and
$\type(w) = \type(\overline{w})$ (we will refer to these as the {\em height} and
{\em type} of $w$).  We wish to determine the relationship between $w$ and
$\type(w)$.  Let $N$ be the normal subgroup of $F_{\{a,t\}}$ generated
by $a$.  The exact sequence
$$1 \longrightarrow N \longrightarrow F_{\{a,t\}} \longrightarrow F_{\{t\}} \longrightarrow 1$$
splits, so we have $F_{\{a,t\}} = N \rtimes F_{\{t\}}$.  This fits into
the following commutative diagram:
\begin{center}
\begin{tabular}{ccccccccc}
$1$ & $\longrightarrow$ & $N$          & $\longrightarrow$ & $N \rtimes F_{\{t\}}$ & $\longrightarrow$ & $F_{\{t\}}$ & $\longrightarrow$ & $1$ \\ 
    &                   & $\downarrow$ &                   & $\downarrow$          &                   & $\parallel$ &                   &     \\
$1$ & $\longrightarrow$ & $\Z^2$       & $\longrightarrow$ & $\Z^2 \rtimes_M \Z$   & $\longrightarrow$ & $\Z$        & $\longrightarrow$ & $1$ \\
\end{tabular}
\end{center}
The map $N \rightarrow \Z^2$ factors through the abelianization $N^{\Abelian}$ of $N$.  Now, it
is well-known (see, e.g., \cite[Exercise 3.2.3-4]{MKS}) that $N$ is the free group on the
generating set
$$\{t^k a t^{-k} : k \in \Z\}$$
The map
$$t^k a t^{-k} \longmapsto z^k$$
therefore defines an isomorphism from $N^{\Abelian}$ to the group $\Z[z,z^{-1}]$ of Laurent polynomials.
Summing up, we have factored the map
$$\type : F_{\{a,t\}} \longrightarrow \Z^2$$
as a composition
$$F_{\{a,t\}} \longrightarrow N \longrightarrow N^{\Abelian} = \Z[z,z^{-1}] \longrightarrow \Z^2$$
Denote by $\unreducedtype(w)$ the image of $w$ in $\Z[z,z^{-1}]$; we will call this the {\em unreduced
type} of $w$.

More concretely, the splitting $F_{\{a,t\}} = N \rtimes F_{\{t\}}$ shows that every
word $w \in F_{\{a,t\}}$ can be expressed as a product
$$w=(\prod\nolimits_{i=1}^n t^{k_i} a^{l_i} t^{-k_i}) t^h$$
with $h,k_i \in \Z$ and $l_i \in \{\pm 1\}$.  Observe that $h(w)=h$.  
Also, $\unreducedtype(w)$ equals the Laurent polynomial
$$\sum\nolimits_{i=1}^n l_i z^{k_i} \in \Z[z,z^{-1}]$$
Since (continuing our systemic confusion of $a$ with the vector $(1,0) \in \Z^2$)
$$\type(t^{k_i} a^{l_i} t^{-k_i}) = l_i M^{k_i} a$$
we have the following lemma:
\begin{lemma}
\label{lemma:unreducedtype}
All $w \in F_{\{a,t\}}$ satisfy $\type(w) = [\unreducedtype(w)(M)] \cdot a$
\end{lemma}

\subsection{Appearance of Types}
We now determine the length of the shortest word with a specified unreduced type
and height.  We begin with some terminology.  Consider an unreduced type
$$t(z)=\sum\nolimits_i c_i z^i \in \Z[z,z^{-1}]$$
with $c_i \in \Z$ and a height $h \in \Z$.  The Laurent polynomial $t(z)$ can
be divided into three different pieces (depending on $h$).  There are two cases.
If $h \geq 0$, we define
\begin{align*}
\Tail_h(t)   & :=\sum\nolimits_{i=-\infty}^{-1} c_i z^i \\
\Center_h(t) & :=\sum\nolimits_{i=0}^h c_i z^i \\
\Head_h(t)   & :=\sum\nolimits_{i=h+1}^\infty c_i z^i \\
\TailLen_h    & :=\Max \{|i| : \textit{$i=0$ or $i<0, c_i \neq 0$}\} \\
\HeadLen_h    & :=\Max \{i-h : \textit{$i=h$ or $i>h, c_i \neq 0$}\}
\end{align*}
If $h \leq 0$, we define
\begin{align*}
\Tail_h(t)   & :=\sum\nolimits_{i=-\infty}^{h-1} c_i z^i \\
\Center_h(t) & :=\sum\nolimits_{i=h}^0 c_i z^i \\
\Head_h(t)   & :=\sum\nolimits_{i=1}^\infty c_i z^i \\
\TailLen_h    & :=\Max \{|i|-|h| : \textit{$i=h$ or $i<h,c_i \neq 0$}\} \\
\HeadLen_h    & :=\Max \{i : \textit{$i=0$ or $i>0,c_i \neq 0$}\}
\end{align*}
We will refer to $\Tail_h(t)$ as the {\em tail}, $\Center_h(t)$ as the {\em center}, and
$\Head_h(t)$ as the {\em head}.  Also, we will call $\TailLen_h(t)$ the {\em length of the tail} and
$\HeadLen_h(t)$ the {\em length of the head}.  Observe that
$$t = \Tail_h(t) + \Center_h(t) + \Head_h(t)$$
Our theorem is the following:

\begin{theorem}
\label{theorem:graysontheorem}
Let $h$ be a height and let
$$t(z)=\sum\nolimits_i c_i z^i \in \Z[z,z^{-1}]$$
be an unreduced type.
Then the shortest word with this unreduced type and height has length
$$2 \TailLen_h(t) + 2 \HeadLen_h(t) + |h| + \sum\nolimits_i |c_i|$$
\end{theorem}
\begin{proof}
We begin by describing an algorithm for determining the unreduced type and height 
of a word $w$ in $\{a^{\pm 1},t^{\pm 1}\}$.  The algorithm keeps track of two pieces of data, the
{\em partial height} $H \in \Z$ and the {\em partial unreduced type} $T \in \Z[z,z^{-1}]$.  
Both are initialized to $0$.  We read $w$ from left to right.  If we read the letter
$t^l$ with $l=\pm 1$, we add $l$ to $H$.  If we read the letter $a^l$ with $l=\pm 1$,
we add $l z^H$ to $T$.  After reading all of $w$, it is clear that $H=h(w)$ and
that $T=\unreducedtype(w)$.

Now consider any word $w$ with the desired height and unreduced type.  Observe that
each $a^{\pm 1}$ in $w$ contributes exactly one term of the form $\pm z^i$.  Hence
$w$ must contain at least $\sum\nolimits_i |c_i|$ letters of the form $a^{\pm 1}$.  To prove
that $w$ is at least as long as the theorem indicates, it is therefore enough
to show that $w$ contains at least $2 \TailLen_h(t) + 2 \HeadLen_h(t) + |h|$ letters
of the form $t^{\pm 1}$.  We first consider the case $h \geq 0$.  In this case, either
$\Tail_h(t)=0$ or $\Tail_h(t)$ must contain a non-zero term of degree $-\TailLen_h(t)$.  This
implies that during our algorithm the partial height $H$ must at some point equal
$-|\Tail_h(t)|$.  Similarly, either $\Head_h(t)=0$ or $\Head_h(t)$ must contain
a non-zero term of degree $h+\HeadLen_h(t)$.  This implies that during our algorithm
the partial height $H$ must at some point equal $h+\HeadLen_h(t)$.  Since $h(w)=h$,
our algorithm must end with $H=h$.  Summing up, the partial height $H$ (which
changes by $\pm 1$ each time a letter of the form $t^{\pm 1}$ is read) starts at $0$,
ends at $h$, at some point equals $-\TailLen_h(t)$, and at some other
point equals $h+\HeadLen_h(t)$.  Clearly at least 
$2 \TailLen_h(t) + 2 \HeadLen_h(t) + |h|$ letters of the form $t^{\pm 1}$ are necessary,
as desired.  The case of $h \leq 0$ is proven in a similar fashion, with the roles
of $\Tail_h(t)$ and $\Head_h(t)$ reversed.

This proves that the indicated expression is a lower bound on the length of a word
with the desired unreduced type and height.  We now prove that this lower bound
is realized.  Like in the proof of the lower bound, the proofs in the cases 
$h \geq 0$ and $h \leq 0$ are similar; we will only consider $h \geq 0$.  In this case,
the following word has the desired length, unreduced type, and height:
$$t^{-\TailLen_h(t)}(\prod\nolimits_{i=-\TailLen_h(t)}^{-1} a^{c_i} t) a^{c_0} (\prod\nolimits_{i=1}^{h+\HeadLen_h(t)} t a^{c_i}) t^{-\HeadLen_h(t)}$$
\end{proof}

\begin{remark}
After proving a version of Theorem \ref{theorem:graysontheorem}, Grayson 
attempts to set up a complicated system of recurrence relations between 
various subsets of the group.  He expresses the growth function as a power 
series whose coefficients are themselves power series.  He demonstrates that 
there is a sort of linear recurrence relation between these (power series) 
coefficients.  He then claims that this is enough to prove that the growth 
series is rational.  However, absent a proof that (say) the first coefficient 
is in fact a rational function this is insufficient.
\end{remark}

Let $T=\Trace(M)$.  Since $M$ is Anosov, it has two distinct real eigenvalues.  
Let $\lambda$ and $\lambda'$ be the eigenvalues with eigenvectors $v$ and 
$v'$.  Let $\alpha, \alpha' \in \R$ be such that 
$$(1,0)=\alpha v + \alpha' v'$$

\begin{theorem}
\label{theorem:equalitytheorem}
Let $w_1$ and $w_2$ be words in $\{a^{\pm 1},t^{\pm 1}\}$.
Then $w_1$ and $w_2$ represent the same element of 
$G$ if and only if $h(w_1)=h(w_2)$ and 
$1-T z+z^2$ divides the Laurent polynomial
$\unreducedtype(w_1)-\unreducedtype(w_2)$.
\end{theorem}
\begin{proof}
Let
\begin{align*}
\unreducedtype(w_1)=&\sum\nolimits_i c_i z^i \\
\unreducedtype(w_2)=&\sum\nolimits_i c'_i z^i
\end{align*}
Observe that with respect to the basis $\{v,v'\}$ Lemma \ref{lemma:unreducedtype} says that we have
\begin{align*}
\reducedtype(w_1)=&(\alpha \sum\nolimits_i c_i \lambda^i, \alpha' 
                    \sum\nolimits_i c_i \lambda'^i) \\
\reducedtype(w_2)=&(\alpha \sum\nolimits_i c'_i \lambda^i, \alpha' 
                    \sum\nolimits_i c'_i \lambda'^i)
\end{align*}
Since $M$ is a $2 \times 2$ matrix with irrational eigenvalues, $\lambda$ and 
$\lambda'$ have the same minimal polynomial as $M$; i.e. $1-Tz+z^2$, whence 
the theorem.
\end{proof}

Consider the set 
$$X=\Z[z, z^{-1}] \times \Z$$
Define a size function on $X$ by setting
$$\|(\sum\nolimits_i c_i z^i,h)\|:=2\TailLen_h(\sum\nolimits_i c_i z^i)+2\HeadLen_h(\sum\nolimits_i c_i z^i)+|h|+1+
                                 \sum\nolimits_i |c_i|$$
Define a partition $P$ on $X$ by the following equivalence relation.
$$(t_1,h_1) \sim (t_2,h_2) \Longleftrightarrow [h_1=h_2 \textit{ and } 1-Tz+z^2 
                                                \textit{ divides } t_1-t_2]$$

\noindent 
We can now state the following important corollary to the above calculations:
\begin{corollary}
\label{corollary:gnearisometricxmodp}
$G$ is near-isometric to $X/P$ (with constant $c=1$).
\end{corollary}

\begin{remark}
The extra $1$ in the definition of the size function on $X$ simplifies the language
$L$ we construct in Section \ref{section:language}, as it forces {\em every} monomial
in the center of a Laurent polynomial to contribute something to the size, even
if it equals $0$.
\end{remark}

\section{The Language}
\label{section:language}
By Corollary \ref{corollary:fellowtravellerrational}, to prove Theorem
\ref{theorem:maintheorem} it is enough to produce a regular language $L$ with a
partition $P'$ satisfying the falsification by fellow traveler property 
so that $L/P'$ is isometric to $X/P$.  We first prove a number of finiteness 
results about $X$.  Next, we will define a series of languages $L_n$ and a series
of corresponding partitions $P_n$.  Finally, we will prove that for $n$ sufficiently
large $L_n/P_n$ is isometric to $X/P$ and satisfies the falsification by fellow traveler
property.

\subsection{Finiteness Lemmas}
The coefficients of the Laurent polynomials associated to elements of $X$
are unbounded.  To apply the theory of finite state automata to $X/P$, we
will first prove a lemma which bounds the coefficients of the Laurent polynomials
associated to elements of minimal size in a single subset in $P$.  We will then prove 
two other lemmas which bound the information we need to keep track of while comparing 
elements of $X$ modulo $P$.

\begin{lemma}
\label{lemma:finiteness1}
Let $x=(t,h) \in X$ be so that
$$\|x\| = \Min \{\|x'\| : \textit{$x=x'$ modulo $P$}\}$$
Then the coefficients $c_i$ of $t$ satisfy $|c_i| < 5|T|$.
\end{lemma}
\begin{proof}
By the definition of $P$, we can for each $i$ add 
or subtract $z^{i-1}-T z^i+z^{i+1}$ from $t$ without changing $\overline{x}$.
Now, if $|c_i| \geq 5|T|$, add or subtract 
$5z^{i-1}-5 T z^i + 5 z^{i+1}$ in such a way as to decrease $|c_i|$.  Examining
the formula for the size of an element of $X$, 
we see that we have subtracted $5|T|$ from the size of $x$ and 
added at most $2+2+5+5=14$.  Since $|T| \geq 3$, we conclude that 
$x$ was not of minimal size, a contradiction.
\end{proof}

\begin{lemma}
\label{lemma:finiteness2}
For every positive integer $A$, there exists some positive integer $B_A$ so that the following 
holds.  For $i=1,2$ let $f_i=\sum\nolimits_j c_{i,j} z^j$ with 
$|c_{i,j}| \leq A$.  Assume that $1-Tz+z^2$ divides $f_1-f_2$.  Then the 
coefficients of $(f_1-f_2)/(1-Tz+z^2)$ are bounded by $B_A$.
\end{lemma}
\begin{proof}
Since $|T| \geq 3$, the largest coefficient which is left when we expand out 
$(1-Tz+z^2)g(z)$ is at least as large as the largest coefficient of $g$.  Hence
we may set $B_A = 2A$
\end{proof}

\begin{lemma}
\label{lemma:finiteness3}
For all positive integers $A$ and $B$, there exists some positive integer
$C_{A,B}$ so that the following holds.  For $i=1,2$ let $f_i=\sum\nolimits_j c_{i,j} z^j$ with
$|c_{i,j}| \leq A$.  Assume that
$$f_1-f_2 = ((1-Tz+z^2)\sum\nolimits_j d_j z^j) + (e_1 z+e_2)$$
with $|d_j| \leq B$.  Then $|e_1|,|e_2| \leq C_{A,B}$.
\end{lemma}
\begin{proof}
Observe that the coefficients of $(1-Tz+z^2)\sum\nolimits_j d_j z^j$ are bounded
by $B(|T|+2)$.  Hence the coefficients of
$$e_1 z+e_2 = f_1-f_2-((1-Tz+z^2)\sum\nolimits_j d_j z^j)$$
are bounded by $C_{A,B}:=2A+B(|T|+2)$.
\end{proof}

\subsection{The Language}
Fix a natural number $n \geq 1$.  Let
$$A_n=\{-n,\ldots,n\} \times \{-1,1,2\}$$
be an alphabet with weighting
$$\phi(c,k)=|c|+|k|$$
Consider the language $L_n$ on $A_n$ whose words are of the following form :
$$(\cdot,2)\cdots(\cdot,2)(\cdot,\pm 1)\cdots
  (\cdot,\pm1)(\cdot,2)\cdots(\cdot,2)$$
We require words $w$ in $L_n$ to satisfy the following conditions.
\begin{enumerate}
\item $w$ must contain at least one letter of the form $(\cdot,\pm 1)$
\item The second entries in all the middle terms of $w$ must be identical.
\item If the common second entry in all the middle terms of $w$ is $-1$, then there must
be at least two such middle terms.
\item If the first or last letters of $w$ equal $(c,2)$, then $c \neq 0$
\end{enumerate}
We also define the language $L_n'$ to
consist of all such words $w$ satisfying conditions 1-3 but not necessarily
4.  Both $L_n$ and $L_n'$ are clearly regular.  Define a map 
$\psi':L_n'\longrightarrow X$ by
\begin{align*}
& \psi'(\prod\nolimits_{i=1}^{n_1}(c_i,2) \prod\nolimits_{i=0}^{n_2}(c_i',1) 
       \prod\nolimits_{i=1}^{n_3}(c_i'',2)) \\
& =(\sum\nolimits_{i=1}^{n_1}c_i z^{i-n_1-1}+\sum\nolimits_{i=0}^{n_2}c_i' z^i + 
    \sum\nolimits_{i=1}^{n_3}c_i'' z^{n_2+i},n_2)
\end{align*}
and
\begin{align*}
& \psi'(\prod\nolimits_{i=1}^{n_1}(c_i,2) \prod\nolimits_{i=0}^{n_2}(c_i',-1) 
       \prod\nolimits_{i=1}^{n_3}(c_i'',2)) \\
& =(\sum\nolimits_{i=1}^{n_1}c_i z^{i-n_1-n_2-1}+\sum\nolimits_{i=0}^{n_2}c_i' z^{i-n_2} + 
    \sum\nolimits_{i=1}^{n_3}c_i'' z^{i},-n_2)
\end{align*}
Let $\psi$ be the restriction of $\psi'$ to $L_n \subset L_n'$.  Observe
that $\psi'$ induces a partition $P_n'$ of $L_n'$ and $\psi$ induces a
partition $P_n$ of $L_n$.  The map $\psi$ is 
clearly a size--preserving map from $L_n$ to $X$, and the fact
that we require that if the sign of the center terms is negative then
there must be at least two center terms forces it to be an injection (this condition
prevents trouble from occurring when $h=0$).  Lemma \ref{lemma:finiteness1} 
implies the following:

\begin{theorem}
\label{theorem:lmodpisometricxmodp}
For $n \geq 5|T|$ the induced map $\overline{\psi}:L_n/P_n \longrightarrow X/P$ is an isometry.
\end{theorem}

We now observe that the tripartite division of words in $L_n'$ reflects
the tail-center-head division of the corresponding Laurent polynomials.
If $w \in L_n'$ with $\psi'(w)=(t,h)$, we define 
\begin{align*}
\Tail(w)   &= \Tail_h(t) \\
\Center(w) &= \Center_h(t) \\
\Head(w)   &= \Head_h(t)
\end{align*}
We will refer to these as the {\em tail}, the {\em center}, and the {\em head} of $w$.
We also define $\TailLen(w)$ and $\HeadLen(w)$ to equal the number
of letters of the form $(c,2)$ at the beginning and end of $w$.  We
remark that if $w$ begins or ends with $(0,2)$, then
$\TailLen(w) \neq \TailLen_h(t)$ or $\HeadLen(w) \neq \HeadLen_h(t)$.

\subsection{The Acceptor}
Fix positive integers $n$ and $i$.  Define a language
\begin{align*}
R_{n,i}=\{ & (w_1,w_2) \in L_n \times L_n : \overline{w}_1=\overline{w}_2 \textit{ in } L_n/P_n \\
           & \textit{ and } |\TailLen(w_1) - \TailLen(w_2)|,|\HeadLen(w_1) - \HeadLen(w_2)| \leq i \}
\end{align*}
This section is devoted to proving the following theorem:
\begin{theorem}
\label{theorem:acceptorregular}
$R_{n,i}$ is a regular language.
\end{theorem}
This has the following immediate corollary:
\begin{corollary}
$R_{n,i}$ is an acceptor for the partition $P_n$ of the language $L_n$.
\end{corollary}

\noindent
To prove Theorem \ref{theorem:acceptorregular}, we need the following lemma:

\begin{lemma}
\label{lemma:acceptorregular}
Define
\begin{align*}
R_n'=\{ & (w_1,w_2) \in L_n' \times L_n' : \overline{w}_1=\overline{w}_2 \textit{ in } L_n'/P_n' \\
           & \textit{ and } \TailLen(w_1) = \TailLen(w_2), \HeadLen(w_1)=\HeadLen(w_2) \}
\end{align*}
Then $R_n'$ is a regular language.
\end{lemma}
\begin{proof}
By Theorem \ref{theorem:regular}, it is enough to construct an
automaton accepting $\Rev(R_n')$, or (to put it in another way) to
construct an automaton which reads $w_1$ and $w_2$ from right to left.
Let $B=B_n$ be the constant from
Lemma \ref{lemma:finiteness2} and $C=C_{n,B_n}$ be the constant from
Lemma \ref{lemma:finiteness3}.  Our strategy will be to imitate
the usual polynomial long division algorithm to divide the difference
between the Laurent polynomials associated to $w_1$ and $w_2$ by
$1-T z+z^2$.  We also will make sure that $w_1$ and $w_2$ ``line up''
properly; that is that they have heads, centers, and tails of the
same length.  

Our automaton has a failure state plus the following
set of states :
\begin{align*}
\{(r,l) : &r=c_1 z + c_2 \textit{ with $c_i \in \Z$ so that $|c_i| \leq C$ and } \\
          &\textit{$l \in \{\Head, \Tail, \Center_1, \Center_{-1}, \Center_{-1,1}\}$}\}
\end{align*}
The second entry in a state keeps track of where we are in $w_1$ and $w_2$ (the
label $\Center_1$ means that the center portion consists of terms of the form $(\cdot,1)$, the
label $\Center_{-1,1}$ means that the center portion consists of terms of the form $(\cdot,-1)$ and we have
only read one term of that form, and the label $\Center_{-1}$ means that the
center portion consists of terms of the form $(\cdot,-1)$ and we have read at least two terms of that
form). The first entry keeps track of the remainder obtained by dividing the
difference of the portion read so far by $1-T z+z^2$.  Recalling that
our automaton reads $w_1$ and $w_2$ from right to left, we begin in
the state $(0,\Head)$.  Assume that $\psi'(w_i) = (t_i,h_i)$ with
$$t_i = \sum\nolimits_{j=-N_1}^{N_2} c_{i,j} z^j$$

Assume now that we are in the state $(r,l)$ after reading $k$ letters.  This means
that there exists some Laurent polynomial $q$ (whose value does not matter -- all that
matters for determining the transitions are the values of $r$ and $l$) so that
$$\sum\nolimits_{i=N_2-k+1}^{N_2} (c_{1,i}-c_{2,i}) z^i = z^{N_2-k+1} (q \cdot (1-T z+z^2) + r)$$
If we do not read entries of the form $(c_{1,N_2-k},e)$ from $w_1$ and $(c_{2,N_2-k},e)$ from
$w_2$ (in other words, if at this point $w_1$ and $w_2$ cease to ``line up''), then we fail.  Otherwise,
the difference between the portions read so far is
\begin{align*}
\sum\nolimits_{i=N_2-k}^{N_2} (c_{1,i}-c_{2,i}) z^i = z^{N_2-k} (&z q \cdot (1-T z+z^2) \\
                                                        &+ (z r + (c_{1,N_2-k}-c_{2,N_2-k})))
\end{align*}
Note that $z r + (c_{1,N_2-k}-c_{2,N_2-k})$ is a quadratic polynomial.  Divide it by $1-T z+z^2$ to get
$$z r + (c_{1,N_2-k}-c_{2,N_2-k}) = q' (1-T z+z^2) + r'$$
where $r'$ is a linear function.  If the coefficients of $r'$ are not bounded by
$C$, then by Lemma \ref{lemma:finiteness3} it is impossible for $w_1$ and $w_2$ to define
equal elements of $L_n'$ modulo $P_n'$, and we fail.  Otherwise, we make the following transition : If
$l=\Head$ or $l=\Tail$ and $e=2$, then we transition to $(r',l)$.  If $l=\Head$ and $e=1$, then
we transition to $(r',\Center_1)$.  If $l=\Head$ and $e=-1$, then we transition to $(r',\Center_{-1,1})$.
If $l=\Center_1$ and $e=1$, we transition to $(r',\Center_1)$.  If $l=\Center_1$ or $l=\Center_{-1}$ and $e=2$, 
we transition to $(r',\Tail)$.  If $l=\Center_{-1,1}$ or $l=\Center_{-1}$ and $e=-1$, we transition to $(r',\Center_{-1})$.  If we
are not in one of these situations, we fail.

Assume now that we manage to successfully read all of $w_1$ and $w_2$ and end in a state $(r,l)$.  
This implies, in particular, that the heads, centers, and tails of $w_1$ and $w_2$ are of the same length.
Also, it is clear from the algorithm that $r$ is the remainder of the difference of the Laurent polynomials
associated to $w_1$ and $w_2$ divided by $1-T z+z^2$.  We succeed if we end in one of the following three states : $(0,\Tail)$,
$(0,\Center_1)$, or $(0,\Center_{-1})$.  The restriction on $l$ is required to guarantee that both $w_1$ and $w_2$ contain
centers of the appropriate form.
\end{proof}

\noindent
We now prove Theorem \ref{theorem:acceptorregular}.

\begin{proofof}{\ref{theorem:acceptorregular}}
Set $Q_r = \prod\nolimits_{j=1}^r (0,2)$.  Observe that
\begin{align*}
R_{n,i} = \{(w_1,w_2) : &\textit{$w_1,w_2 \in L_n$ and there exist $r,s \in \Z$ so that $0 \leq r,s \leq i$} \\
                        &\textit{and either $(Q_r w_1 Q_s, w_2) \in R_n'$, $(Q_r w_1, w_2 Q_s) \in R_n'$,} \\
                        &\textit{$(w_1 Q_s, Q_r w_2) \in R_n'$, or $(w_1,Q_r w_2 Q_s) \in R_n'$} \}
\end{align*}
Since the integers $r$ and $s$ which appear in this expression are bounded, it can be
expressed using first order predicates and concatenation.  Theorem \ref{theorem:regular} therefore 
implies that $R_{n,i}$ is a regular language.
\end{proofof}

\subsection{The Falsification by Fellow Traveler Property and Proof of the Main Theorem}
In this section we will complete the proof of Theorem \ref{theorem:maintheorem}.  We
will need the following definition:

\begin{definition}
Consider $w_1,w_2 \in L_n$ with $\psi(w_i)=(\sum\nolimits_j c_{i,j} z^j,h_i)$.  The {\em divergence} of
$w_1$ and $w_2$ is the maximal absolute value of
$$\sum\nolimits_{j=-\infty}^k |c_{1,j}|-|c_{2,j}|$$
as $k$ varies.
\end{definition}

\noindent
The key step in our proof will be the following theorem:

\begin{theorem}
\label{theorem:weaklyfellowtravel}
There exist constants $K$, $L$, and $N$ so that $N \geq 5|T|$ and the following are true.
\begin{enumerate}
\item If $w_1 \in L_{5|T|}$ is not a minimal size representative modulo $P_{5|T|}$, then
there exists some $w_2 \in L_N$ so that 
\begin{itemize}
\item $\overline{w}_1=\overline{w}_2$ and $\|w_2\| < \|w_1\|$.
\item $|\HeadLen(w_1) - \HeadLen(w_2)| \leq L$ and $|\TailLen(w_1) - \TailLen(w_2)| \leq L$
\item The divergence of $w_1$ and $w_2$ is bounded by $K$.
\end{itemize}
\item If $w_1,w_2 \in L_{5|T|}$ are two different minimal size representatives of the same
element modulo $P_{5|T|}$, then $|\HeadLen(w_1) - \HeadLen(w_2)| \leq L$ and $|\TailLen(w_1) - \TailLen(w_2)| \leq L$.
\end{enumerate}
\end{theorem}

\noindent
Before proving Theorem \ref{theorem:weaklyfellowtravel}, we will use it to prove Theorem \ref{theorem:maintheorem}.

\begin{proofof}{\ref{theorem:maintheorem}}
By Corollary \ref{corollary:gnearisometricxmodp} and Theorem \ref{theorem:lmodpisometricxmodp}, it is
is enough to show that $L_{N} / P_{N}$ has a rational growth series for large $N$.  Let $K$, $L$, and $N$ be 
the constants from Theorem \ref{theorem:weaklyfellowtravel}.  We will prove
that $P_N$ is a partition of $L_N$ with acceptor $R_{N,L}$ satisfying the falsification by fellow traveler
condition with respect to the constant $K+(N+6)L$.  By Theorem \ref{corollary:fellowtravellerrational}, this
will imply that $L_N / P_N$ has a rational growth series, as desired.  

We begin by observing that by Lemma \ref{lemma:finiteness1}, $L_{5|T|}$ contains minimal size elements from 
each set in $P_N$.  Now, let $w_1 \in L_{5|T|}$ not be
a minimal size representative modulo $P_{5|T|}$.  Consider the $w_2 \in L_N$ given by
the first conclusion of Theorem \ref{theorem:weaklyfellowtravel}.  It is clear that $(w_1,w_2) \in R_{N,L}$
and that $\|w_2\| < \|w_1\|$.  We must prove that $w_1$ and $w_2$ $(K+(N+6)L)$-fellow travel.  

We will assume that $\TailLen(w_2) \leq \TailLen(w_1)$; the other case is similar.  Consider 
length $j$ initial segments $v_1$ and
$v_2$ of $w_1$ and $w_2$.  Let $v_2'$ be the initial segment of $w_2$ of length 
$j-(\TailLen(w_1)-\TailLen(w_2))$.  Since the divergence of $w_1$ and $w_2$ is bounded by $K$, we know that
$$| \|v_1\| - \|v_2'\| | \leq K + 4L$$
The $4L$ term comes from the fact that each term in the initial segment of length $\TailLen(w_1)-\TailLen(w_2) \leq L$
of $v_1$ contributes an extra $2$ to the difference, and in addition either $v_1$ or $v_2'$ may contain
at most $L$ terms from the head which are absent from the other, each possibly contributing $2$ more to the
difference.  
The remaining portion of $v_2$ has length at most $L$ and each term
contributes at most $2+N$ to the size of $v_2$.  Hence we conclude that
$$| \|v_1\| - \|v_2\| | \leq K+4L + L(2+N) = K+(N+6)L$$
as desired.

Now, by the second conclusion of Theorem \ref{theorem:weaklyfellowtravel}, if $w_1$ and $w_2$ are 
two minimal-size representatives of the same element of $L_{5|T|}$ modulo $P_{5|T|}$, then $(w_1,w_2) \in R_{N,L}$.  This completes the
proof of the falsification by fellow traveler property, and hence of the theorem.
\end{proofof}

\noindent
We now prove Theorem \ref{theorem:weaklyfellowtravel}.

\begin{proofof}{\ref{theorem:weaklyfellowtravel}}
Let $B:=B_{5|T|}$ be the constant from Lemma \ref{lemma:finiteness2}.  We will
prove that the following choices of $L$, $K$, and $N$ suffice:
\begin{align*}
L&=(|T|+2)B \\
K&=(|T|+2)(3B+4)+8L+1 \\
N&=5|T|+(|T|+2)B
\end{align*}
We begin by proving the first conclusion of the theorem.  Let $w_1 \in L_{5|T|}$ not be
a minimal size representative modulo $P_{5|T|}$.  By Lemma \ref{lemma:finiteness1},
there exists some $w_2 \in L_{5|T|}$ so that $\overline{w}_1=\overline{w}_2$ and
$\|w_2\| < \|w_1\|$.  Let $t_1,t_2 \in \Z[z,z^{-1}]$ and $h \in \Z$ be
so that $(t_i,h) = \psi(w_i)$.  Expand the $t_i$ as $t_i = \sum\nolimits_j c_{i,j} z^j$.
By Lemma \ref{lemma:finiteness2}, there exists some Laurent polynomial $q=\sum\nolimits_j d_j z^j$
with $|d_j| \leq B$ so that $t_2 = t_1 + (1-T z+z^2) q$.

Our goal will be to modify $t_2$ and $q$ to produce new Laurent polynomials
$t_2'$ and $q'$ with $t_2' = t_1 + (1-T z+z^2) q'$ so that
(expanding $q'$ and $t_2'$ as $q'=\sum\nolimits_j d_j' z^j$ and $t_2'=\sum\nolimits_j c_{2,j}' z^j$ and 
setting $w_2'=\psi^{-1}(t_2',h)$) the following conditions are satisfied:
\begin{enumerate}
\item $|d_j'| \leq B$ and $\|w_2'\| < \|w_1\|$
\item $\TailLen(w_2') - \TailLen(w_1) \leq L$ and $\HeadLen(w_2') - \HeadLen(w_1) \leq L$
\item $\TailLen(w_1) - \TailLen(w_2') \leq L$ and $\HeadLen(w_1) - \HeadLen(w_2') \leq L$
\item For all $k$ we have $\sum\nolimits_{j=-\infty}^k (|c_{2,j}'| - |c_{1,j}|) \leq K$.
\item For all $k$ we have $\sum\nolimits_{j=-\infty}^k (|c_{1,j}| - |c_{2,j}'|) \leq K$.
\end{enumerate}
The first part of condition $1$ implies that $w_2' \in L_N$, and the rest of the conditions
imply the first conclusion of the theorem.  Our modification will take several steps; to prevent
a proliferation of new notation we will continue to refer to the modified polynomials, words, and
coefficients as $t_2$, $q$, $w_2$, $d_j$, and $c_{2,j}$.  The modifications are the following:

\begin{claim}{1}
We can modify $w_2$ so that conditions 1 and 2 are satisfied.
\end{claim}

\begin{claimproof}
We will indicate how to achieve $\TailLen(w_2) \leq \TailLen(w_1) + L$; the other
modification is similar.  Assume that $\TailLen(w_2) > \TailLen(w_1) + (|T|+2)B$.  We will
show how to find a $w_2'$ so that condition 1 is satisfied and so that 
$\|w_2'\| < \|w_2\|$; repeating this process will eventually yield the desired conclusion.  The
idea of our construction is that since each element of the tail of $w_2$ contributes something
to $\|w_2\|$, if the tail is sufficiently long then we can remove the first few terms from it and
shrink $\|w_2\|$.  Let $M$ be the smallest integer with $d_M \neq 0$.  Hence
$$t_2 = t_1 + (1-T z+z^2) \sum\nolimits_{j=M}^{\infty} d_j z^i$$
Set $M' = M+(|T|+2)B$ and
$$t_2' = t_1 + (1-T z+z^2)\sum\nolimits_{j=M'}^{\infty} d_j z^i$$
Expand this as $t_2'=\sum\nolimits_j c_{2,j}' z^j$ and set $w_2'=\psi^{-1}(t_2',h)$.

The only non-trivial fact we must prove is $\|w_2'\| < \|w_2\|$.  
Observe first that by construction we have
$$\TailLen(w_2) - \TailLen(w_2') \geq (|T|+2)B$$
Also, 
$$|c_{2,j}'| =
\begin{cases}
0                                                    & \text{for $j < M'$} \\
|c_{2,j}-d_{j-2}+T d_{j-1}| \leq |c_{2,j}|+(|T|+1)B & \text{for $j=M'$} \\
|c_{2,j}-d_{j-2}| \leq |c_{2,j}|+B        & \text{for $j=M'+1$} \\
|c_{2,j}|                                            & \text{for $j > M'+1$}
\end{cases}$$
Finally, we may have lengthened the head of $w_2$ by 1; i.e.
$$\HeadLen(w_2) - \HeadLen(w_2') \geq -1$$
Summing up,
\begin{align*}
\|w_2\| - \|w_2'\| &=    2(\TailLen(w_2) - \TailLen(w_2')) + 2(\HeadLen(w_2)-\HeadLen(w_2')) + \sum\nolimits_j |c_{2,j}|-|c_{2,j}'|  \\
                   &\geq 2(|T|+2)B - 2 - (|T|+1)B - B = (|T|+2)B - 2 > 0
\end{align*}
as desired.
\end{claimproof}

\begin{claim}{2}
We can modify the $w_2$ produced in Claim 1 so that conditions 1-3 are satisfied.
\end{claim}

\begin{claimproof}
Assume that $\TailLen(w_1)-\TailLen(w_2) > L$.  The idea of our construction of $w_2'$ is that
since each term of the tail of $w_1$ contributes something to $\|w_1\|$ and $w_2$ has a much
shorter tail than $w_1$ we can use only the initial portion of the quotient $q$ to shorten the tail of
$w_1$ by enough to shrink $\|w_1\|$.  Let $M$ be the smallest integer 
with $c_{1,M} \neq 0$.  Set $M'=M+(|T|+2)B$ and let
$$t_2' = t_1 + (1-T z+z^2)\sum\nolimits_{j=-\infty}^{M'-1} d_j z^i$$
Expand this as $t_2'=\sum\nolimits_j c_{2,j}' z^j$ and set $w_2'=\psi^{-1}(t_2',h)$.

Observe that $c_{2,j}'=0$ for $j<M'$.  Informally, we have ``chopped off'' 
the first $(|T|+2)B$ terms from the tail of $w_1$.  However,
we may have been too successful : possibly $c_{2,M'}' = 0$, indicating that we have
shortened the tail more than we intended.  If this is the case, add or subtract $z^{M'}(1-Tz+z^2)$
from $t_2'$ in such a way as to insure that we still have $w_2' \in L_N$.  There is
therefore some integer $E \in \{-1,0,1\}$ so that
\begin{align*}
|c_{2,j}'| &= 
\begin{cases}
0                               & \text{if $j < M'$} \\
|c_{1,j}+d_{j-2}-T d_{j-1} + E| & \text{if $j=M'$} \\
|c_{1,j}+d_{j-2} - E T|         & \text{if $j=M'+1$} \\
|c_{1,j}+E|                     & \text{if $j=M'+2$} \\
|c_{1,j}|                       & \text{if $j > M'+2$}
\end{cases} \\
&\leq
\begin{cases}
0                    & \text{if $j < M'$} \\
|c_{1,j}|+(|T|+1)B+1 & \text{if $j=M'$} \\
|c_{1,j}|+B+|T|      & \text{if $j=M'+1$} \\
|c_{1,j}|+1          & \text{if $j=M'+2$} \\
|c_{1,j}|            & \text{if $j > M'+2$}
\end{cases}
\end{align*}
Also,
$$\TailLen(w_1)-\TailLen(w_2') = (|T|+2)B = L$$
Finally, we may have changed the length of the head of $w_1$ by 1; i.e.
$$|\HeadLen(w_1) - \HeadLen(w_2')| \leq 1$$
These facts imply that $w_2'$ satisfies condition 2 and 3.  To show that $w_2'$ also satisfies
condition 1, we must show that $\|w_2'\|<\|w_1\|$.  This follows from the 
following calculation.
\begin{align*}
\|w_1\| - \|w_2'\| &=    2(\TailLen(w_1) - \TailLen(w_2')) + 2(\HeadLen(w_1)-\HeadLen(w_2')) + \sum\nolimits_j |c_{1,j}|-|c_{2,j}'|  \\
                   &\geq 2(|T|+2)B - 2 - ((|T|+1)B+1) - (B+|T|) - 1\\
                   &=    (|T|+2)B-|T|-4 = |T|(B-1) + (2B-4) > 0
\end{align*}
The final inequality follows from the fact that $B \geq 2$.  In a similar way, one can show that if $\HeadLen(w_1) - \HeadLen(w_2) > L$ then one
can modify $w_2$ in an appropriate way.
\end{claimproof}

\begin{claim}{3}
We can modify the $w_2$ produced in Claim 2 so that conditions 1-4 are satisfied.
\end{claim}

\begin{claimproof}
Assume that for some $k$ we have
$$\sum\nolimits_{j=-\infty}^k (|c_{2,j}| - |c_{1,j}|) > (|T|+2)B+8L+1 + 2(|T|+2)$$
We will show that we can find some $w_2''$ satisfying conditions 1-3 so that
$\|w_2''\| < \|w_2\|$; repeating this process will eventually yield the desired 
conclusion.  This construction will be a two-step process.  The idea of the first
part of our construction is that since the initial portion
of $w_2$ is so much larger than the corresponding portion of $w_1$ we can remove
the initial portion from $q$ to get a $w_2'$ which begins like $w_1$ and ends
like $w_2$ and is smaller than $w_2$.  Set
$$t_2' = t_1 + (1-T z+z^2)\sum\nolimits_{j=k+1}^\infty d_j z^j$$
Expand this as $t_2' = \sum\nolimits_j c_{2,j}' z^j$
and set $w_2' = \psi^{-1}(t_2',h)$.  Observe that
$$|c_{2,j}'| =
\begin{cases}
|c_{1,j}|                                                 & \text{if $j \leq k$} \\
|c_{2,j} - d_{j-2} + T d_{j-1}| \leq |c_{2,j}| + (|T|+1)B & \text{if $j=k+1$} \\
|c_{2,j} - d_{j-2}| \leq |c_{2,j}| + B                    & \text{if $j=k+2$} \\
|c_{2,j}|                                                 & \text{if $j > k+2$}
\end{cases}$$
Now, in a manner similar to that in Claim 2, we may have inadvertently shortened
the tail or head of $w_2$ so much that $w_2'$ no longer satisfies condition 3.  To fix this,
create a new Laurent
polynomial $t_2''$ by adding $(E_1 z^{M_1} + E_2 z^{M_2})(1-Tz+z^2)$ with
$E_1,E_2 \in \{-1,0,1\}$ and $M_1,M_2 \in \Z$ to $t_2'$ in such a way as to assure that (setting
$w_2'' = \psi^{-1}(t_2'',h)$) we have $w_2'' \in L_N$ and
$$|\HeadLen(w_1) - \HeadLen(w_2'')|, |\TailLen(w_1) - \TailLen(w_2'')| \leq L$$
Observe that the divergence of $w_2'$ and $w_2''$ is bounded by $2(|T|+2)$.  
This implies that
\begin{align*}
\|w_2\| - \|w_2''\| \geq &2(\TailLen(w_2)-\TailLen(w_2'')) + 2(\HeadLen(w_2)-\HeadLen(w_2'')) \\
                         &+ (\sum\nolimits_j (|c_{2,j}|-|c_{2,j}'|) - 2(|T|+2)) \\
\geq &2(\TailLen(w_2)-\TailLen(w_1)) + 2(\TailLen(w_1) - \TailLen(w_2'')) \\
     &+ 2(\HeadLen(w_2)-\HeadLen(w_1)) + 2(\HeadLen(w_1) - \HeadLen(w_2'')) \\
     &+\sum\nolimits_{j=-\infty}^k (|c_{2,j}| - |c_{1,j}|) - (|T|+1)B - B - 2(|T|+2)\\
>    &-8L + (|T|+2)B+8L+1+2(|T|+2) \\
     &- (|T|+1)B - B - 2(|T|+2)\\
=    &1
\end{align*}
as desired.
\end{claimproof}

\begin{claim}{4}
We can modify the $w_2$ produced in Claim 3 so that conditions 1-5 are satisfied.
\end{claim}

\begin{claimproof}
Assume that for some $k$ we have
$$\sum\nolimits_{j=-\infty}^k (|c_{1,j}| - |c_{2,j}|) > (|T|+2)B+8L+1+2(|T|+2)$$
Pick this $k$ to be the minimal $k$ with this property.  Since                       
$|c_{1,j}| - |c_{2,j}| \leq (|T|+2)B$, we have
$$\sum\nolimits_{j=-\infty}^k |c_{1,j}| - |c_{2,j}| < ((|T|+2)B+8L+1)+2(|T|+2)+(|T|+2)B$$
We will construct from $w_2$ a $w_2''$ satisfying conditions 1-4 whose divergence from $w_1$
is bounded by $K$.  Again, the construction of $w_2''$ is a two-step process.  The idea of the first part of
our construction is that since the initial portion of $w_2$ is so much
smaller than the corresponding portion of $w_1$ we can use only the initial portion of $q$ to get a word $w_2'$ which
is definitely smaller than $w_1$.  Set
$$t_2' = t_1 + (1-T z+z^2)\sum\nolimits_{j=-\infty}^k d_j z^j$$
Expand this as $t_2' = \sum\nolimits_j c_{2,j}' z^j$ and set $w_2'=\psi^{-1}(t_2',h)$.  Observe that
$$c_{2,j}' =
\begin{cases}
c_{2,j} & \text{if $j \leq k$} \\
c_{1,j} & \text{if $j > k+2$}
\end{cases}$$
and
$$||c_{2,j}'| - |c_{1,j|}|| \leq
\begin{cases}
(|T|+1)B & \text{if $j=k+1$} \\
B        & \text{if $j=k+2$}
\end{cases}$$
Now, like in Claim 3 we may have inadvertently shortened
the tail or head of $w_2$ so much that $w_2'$ no longer satisfies condition 3.  To fix this,
create a new Laurent
polynomial $t_2''$ by adding $(E_1 z^{M_1} + E_2 z^{M_2})(1-Tz+z^2)$ with
$E_1,E_2 \in \{-1,0,1\}$ and $M_1,M_2 \in \Z$ to $t_2'$ in such a way as to assure that (setting
$w_2'' = \psi^{-1}(t_2'',h)$) we have $w_2'' \in L_N$ and
$$|\HeadLen(w_1) - \HeadLen(w_2'')|, |\TailLen(w_1) - \TailLen(w_2'')| \leq L$$
Observe that the divergence of $w_2'$ and $w_2''$ is bounded by $2(|T|+2)$.
This implies that the divergence of $w_1$ and $w_2''$ is bounded by $2(|T|+2)$ plus the divergence
of $w_1$ and $w_2''$.  The above formulas plus the minimality of $k$ imply that the
divergence of $w_1$ and $w_2'$ is bounded by
\begin{align*}
  &((|T|+2)B+8L+1)+2(|T|+2)+(|T|+2)B+(|T|+1)B+B \\
= &(|T|+2)(3B+2)+8L+1
\end{align*}
We conclude that the divergence of $w_1$ and $w_2''$ is bounded by
$$(|T|+2)(3B+2)+8L+1+2(|T|+2) = (|T|+2)(3B+4)+8L+1 = K$$
as desired.  It is enough, therefore, to prove that $\|w_2''| < \|w_1\|$.  By the above,
\begin{align*}
\|w_1\|-\|w_2''\| \geq &2(\TailLen(w_1) - \TailLen(w_2')) + 2(\HeadLen(w_1)-\HeadLen(w_2')) \\
                       &+ (\sum\nolimits_j (|c_{1,j}|-|c_{2,j}'|) - 2(|T|+2)  \\
                  >    &-2L-2L+((|T|+2)B+8L+1+2(|T|+2) \\
                       &- 2(|T|+2))-(|T|+1)B-B \\
                  =    &4L+1>0
\end{align*}
as desired.
\end{claimproof}

\noindent
These claims complete the proof of the first conclusion of the theorem.  

We now prove the second conclusion.  We recall that the second conclusion of the theorem
is that if $w_1$ and $w_2$ are two minimal size representatives of the same
element modulo $P_{5|T|}$, then $|\HeadLen(w_1) - \HeadLen(w_2)| \leq L$ and $|\TailLen(w_1) - \TailLen(w_2)| \leq L$.
Assume that $w_1$ and $w_2$ are equal modulo $P_{5|T|}$ and satisfy either 
$|\HeadLen(w_1) - \HeadLen(w_2)| > L$ or $|\TailLen(w_1) - \TailLen(w_2)| > L$.  Without loss
of generality assume that either $\HeadLen(w_2) - \HeadLen(w_1) > L$ or $\TailLen(w_2) - \TailLen(w_1) > L$.
The proof of Claim 1 tells us then that there exists some $w_2'$ satisfying $\overline{w}_2'=\overline{w}_1$
and $\|w_2'\| < \|w_2\|$.  In particular, $w_2$ was not of minimal size, as desired.
\end{proofof}

\section{Some Questions}
\label{section:questions}
As we remarked in the introduction, using the methods of this paper to actually
compute growth series would be a long and unpleasant task.  However, in an 
unpublished paper Parry has calculated some growth series for torus bundle 
groups \cite{P}.  We reproduce his formulas here.  He considers a torus bundle 
whose monodromy has trace $2T$.  Letting $\langle a,b,t \rangle$ be the natural
generators, he calculates the growth series of the finite--index subgroup
generated by $\langle a,t a t^{-1},t \rangle$ with respect to that generating 
set (observe that this is the same subgroup we considered, but with one 
additional generator.  It is not too hard to adapt our proof to this new 
generating set).  He proves that the growth function is $N(z)/D(z)$, where 
$N(z)$ and $D(z)$ are the following.
\begin{align*}
N(z)=&(1-z)^2 (1+z) (1+3z+4z^2+4z^3+3z^4+z^5 \\
     &-z^T-3z^{T+1}-14z^{T+2}-16z^{T+3}-11z^{T+4}-5z^{T+5}+2z^{T+6} \\
     &+2z^{2T+1}-13z^{2T+2}+35z^{2T+3}+40z^{2T+4}+6z^{2T+5} \\
     &-23z^{2T+6}-7z^{2T+7}+4z^{2T+8}+4z^{2T+9} \\
     &-5z^{3T+2}+31z^{3T+3}-40z^{3T+4}-44z^{3T+5} \\
     &+33z^{3T+6}+25z^{3T+7}-12z^{3T+8}-4z^{3T+9}) \\
D(z)=&(1-2z-z^2-z^T+4z^{T+1}-z^{T+2})(1-z-z^2-z^3 \\
     &-z^{T+1}+3z^{T+2}+z^{T+3}-z^{T+4})^2
\end{align*}
Since Parry only dealt with the even trace case, we pose the following 
combinatorial challenge.

\smallskip
\noindent
{\bf Question 1:} Explicitly compute the growth series of our finite index 
subgroups for torus bundles with odd trace monodromy.
\smallskip

The 3--dimensional Sol groups are the fundamental groups of 2 dimensional
torus bundles over a circle whose monodromy has has no eigenvalues on the
unit circle.  By considering $n$ dimensional torus bundles over a circle with
the same restriction on the monodromy, we get the $n+1$--dimensional Sol groups.
It seems difficult to generalize our methods to these groups.  This suggests
the following question.

\smallskip
\noindent
{\bf Question 2:} Do the higher--dimensional Sol groups have rational growth
functions?
\smallskip

The fact that we were only able to find a finite index subgroup with rational 
growth suggests the following question.

\smallskip
\noindent
{\bf Question 3:} Does there exist any group $G$ which has irrational growth 
with respect to all sets of generators but which contains a finite index 
subgroup $G'$ which has rational growth with respect to some set of generators?
\smallskip

\noindent 
The following more general question also seems interesting.

\smallskip
\noindent
{\bf Question 4:} Consider the property of having rational growth with respect 
to some set of generators.  How does this property behave under 
commensuration?  under quasi-isometry?
\smallskip

\begin{remark}
Observe that $\widetilde{\SL_2}$ and $\HBolic^2 \times \R$ are quasi-isometric,
and hence the fundamental groups of manifolds modeled on these geometries are 
quasi-isometric.  An easy consequence of Cannon's work (see \cite{Ca}) is that 
the fundamental groups of manifolds modeled on $\HBolic^2 \times \R$ are 
rational with respect to any generating set.  The work of Shapiro suggests that
this is probably false for manifolds modeled on $\widetilde{\SL_2}$ (see 
\cite{Sh2}), so the property of being rational with respect to {\em all\/} 
generating sets is likely not well--behaved under quasi-isometry.  The 
question of whether the fundamental group of any $\widetilde{\SL_2}$--manifold 
has rational growth with respect to some generating set is still open.
\end{remark}

\noindent
Dept. of Mathematics, University of Chicago\\
5734 University Ave.\\
Chicago, Il 60637\\
E-mail: {\tt andyp@math.uchicago.edu}
\medskip

\end{document}